\documentclass[12pt,fleqn]{article}
\usepackage{graphicx}

\usepackage{latexsym}

\usepackage{amsmath}
\usepackage{amsthm}
\usepackage{amssymb}
\usepackage{amsfonts}

\numberwithin{equation}{section}

\begin{document}

\newcommand{\rf}[1]{(\ref{#1})}
\newcommand{\rff}[2]{(\ref{#1}\ref{#2})}

\newcommand{\ba}{\begin{array}}
\newcommand{\ea}{\end{array}}

\newcommand{\be}{\begin{equation}}
\newcommand{\ee}{\end{equation}}

\newcommand{\const}{{\rm const}}
\newcommand{\ep}{\varepsilon}
\newcommand{\Cl}{{\cal C}}
\newcommand{\rr}{{\vec r}}
\newcommand{\ph}{\varphi}
\newcommand{\R}{{\mathbb R}}  
\newcommand{\N}{{\mathbb N}}
\newcommand{\Z}{{\mathbb Z}}

\newcommand{\e}{{\bf e}}

\newcommand{\m}{\left( \ba{r}}
\newcommand{\ema}{\ea \right)}
\newcommand{\mm}{\left( \ba{cc}}
\newcommand{\miv}{\left( \ba{cccc}}

\newcommand{\scal}[2]{\mbox{$\langle #1 \! \mid #2 \rangle $}}
\newcommand{\ods}{\par \vspace{0.3cm} \par}
\newcommand{\dis}{\displaystyle }
\newcommand{\mc}{\multicolumn}
\newcommand{\no}{\par \noindent}

\newcommand{\sinc}{ {\rm sinc\,} }
\newcommand{\tanhc}{{\rm tanhc} }
\newcommand{\tanc}{{\rm tanc} }

\newtheorem{prop}{Proposition}[section]
\newtheorem{Th}[prop]{Theorem}
\newtheorem{lem}[prop]{Lemma}
\newtheorem{rem}[prop]{Remark}
\newtheorem{cor}[prop]{Corollary}
\newtheorem{Def}[prop]{Definition}
\newtheorem{open}{Open problem}
\newtheorem{ex}[prop]{Example}
\newtheorem{exer}[prop]{Exercise}

\newenvironment{Proof}{\par \vspace{2ex} \par
\noindent \small {\it Proof:}}{\hfill $\Box$ 
\vspace{2ex} \par }

\title{\bf 
Locally exact modifications of numerical integrators }
\author{
 {\bf Jan L.\ Cie\'sli\'nski}\thanks{\footnotesize
 e-mail: \tt janek\,@\,alpha.uwb.edu.pl}
\\ {\footnotesize Uniwersytet w Bia{\l}ymstoku,
Wydzia{\l} Fizyki}
\\ {\footnotesize ul.\ Lipowa 41, 15-424
Bia{\l}ystok, Poland}
}

\date{}

\maketitle

\begin{abstract}
We present a new class of exponential integrators for ordinary differential equations. They are locally exact, i.e., they preserve the linearization of the original system at every point. Their construction consists in modifying existing numerical schemes in order to make them locally exact. The resulting schemes preserve all fixed points and are A-stable. The most promising results concern the discrete gradient method (modified implicit midpoint rule) where we succeeded to preserve essential geometric properties and the final results have a relatively simple form. In the case of one-dimensional Hamiltonian systems numerical experiments show that our modifications can increase the accuracy by several orders of magnitude. The main result of this paper is the construction of energy-preserving locally exact discrete gradient schemes for arbitrary multidimensional Hamiltonian systems in canonical coordinates. 
\end{abstract}

\ods




\section{Introduction}

The motivation for introducing ``locally exact'' discretizations is quite natural. Considering a numerical scheme for a dynamical system with periodic solutions (for instance: a nonlinear pendulum) we may ask whether the numerical scheme recovers the period of small oscillations. For a fixed finite time step the answer is usualy negative. However, many numerical scheme (perhaps all of them) admit modifications which preserve the period of small oscillations for a fixed (not necessarily small) time step $h$.  
An unusual feature of our approach is that instead of taking the limit $h \rightarrow 0$, we consider the limit $x_n \approx \bar x$ (where $\bar x$ is the stable equilibrium). 
The next step is to consider a linearization around any fixed $\bar x$ (then, in order to preserve the condition $x_n \approx \bar x$, some evolution of $\bar x$ is necessary). 
In other words, we combine two known procedures: the approximation of nonlinear systems by linear equations and explicit exact discretizations of linear equations with constant coefficients. An important novelty consists in applying these procedures to a modified numerical scheme containing free functional parameters.  The case of small oscillations was presented in \cite{CR-long}. More results for one-dimensional Hamiltonian systems can be found in \cite{CR-PRE,CR-BIT} (by one-dimensional Hamiltonian system we mean a Hamiltonian system with one degree of freedom). The results are very promising. It seems that a new, very accurate method is emerging. In this paper we extend our approach on the case of multidimensional canonical Hamiltonian systems. 
Our method works perfectly for discrete gradient schemes. We succeeded  to modify the discrete gradient scheme in a locally exact way without spoiling its main geometric property: the exact conservation of the energy integral. It is well known that preservation of geometric properties by numerical algorithms is of considerable advantage \cite{HLW,Is}. 

A notion identical with our local exactness has been proposed a long time ago \cite{Pope}, see also \cite{Law}.  Recently, similar concept appeared under the name of linearization-preserving preprocessing \cite{MQTse}, see also  below (section~\ref{sec-locex}). 
Our approach has  also some similarities with the Mickens approach \cite{Mic}, Gautschi-type methods \cite{Gau,HoL2} and, most of all, with the exponential integrators technique \cite{CCO,HoL1,MW}. The definition of exponential integrators is so wide (e.g., ``a numerical method which involves an exponential function of the Jacobian'', \cite{HoL1}) that our numerical schemes can be considered as special exponential integrators. 
  In spite of some similarities and overlappings,  our approach differs from all methods mentioned above. In particular, according to our best knowledge, discrete gradient schemes have never been treated or modified in the framework of exponential integrators and/or linearization-preserving preprocessing.  Main new results of this paper consist in constructing energy-preserving locally exact discrete gradient schemes for arbitrary multidimensional Hamiltonian systems, see section~\ref{sec-locex-multi}.

\section{Exact discretization of linear systems}
\label{sec-exact}

We consider an ordinary differential equation (ODE) with the general solution ${\pmb x} (t)$ (satisfying the initial condition ${\pmb x} (t_0) = {\pmb x}_0$), and a difference equation with the  general solution ${\pmb x}_n$. 
The difference equation 
 is the exact discretization of  the considered ODE  if \ ${\pmb x}_n = {\pmb x} (t_n)$.

It is well known that any linear ODE with constant coefficients admits the exact discretization in an explicit form \cite{Potts}, see also \cite{Ag,CR-ade,Mic}. We summarize these results as follows, compare \cite{CR-BIT} (Theorem 3.1). 

\begin{prop} \label{prop-exact} 
Any linear equation with constant coefficients, represented in the matrix form by
\be  \label{genlin}
  \frac{d  {\pmb x}}{d t} = A {\pmb x} + {\pmb b} \ ,
\ee
(where ${\pmb x} = {\pmb x} (t) \in \R^d$, ${\pmb b} = \const \in \R^d$ and $A$ is a constant invertible real $d\times d$ matrix) admits the exact discretization given by
\be   \label{genex}
   {\pmb x}_{n+1} - {\pmb x}_n  = (e^{h_n A} - 1) A^{-1} \left( A {\pmb x}_n + {\pmb b} \right)   \ , 
\ee
where \ $h_n = t_{n+1} - t_n$ \  is the time step and $1$ is the identity matrix. 
\end{prop}

\begin{cor}  We may look at the exact discretization \rf{Del-ex} as a modification of the forward Euler scheme. Indeed,  we can rewrite 
 \rf{genex} as 
\be  \label{Del-ex}
  \Delta_n^{-1} \left(  {\pmb x}_{n+1} - {\pmb x}_n  \right) = A {\pmb x}_n + {\pmb b} \ ,
\ee
where  $\Delta_n$, defined by  
\be  
   \Delta_n =  A^{-1} (e^{h_n A} - 1) \ , 
\ee
is  a matrix ``perturbation'' of the time step $h_n$. Indeed,  $ \Delta_n = h_n 1 + O (h_n^2)$. 
\end{cor}

As a special case we consider the multidimensional harmonic oscillator equation driven by a constant force:
\be  \label{vec-osc}
 \frac{d^2   {\pmb x } }{d t^2}  +  \Omega^2 {\pmb x} = {\pmb a} \ , 
\ee
where ${\pmb x} = {\pmb x}(t) \in \R^m$ and $ \Omega$ is a given invertible constant $m \times m$ matrix and ${\pmb a} = \const \in \R^m$. It is convenient to represent 
\rf{vec-osc} as the following first order system
\be  \label{vec1-osc}
 {\pmb {\dot x}} = {\pmb  p} \ , \quad {\pmb {\dot p}} = - \Omega^2 {\pmb x} + {\pmb  a} \ . 
\ee
 
\begin{prop}[\cite{Ci-oscyl}, Prop.\ 19]
The exact discretization of the system \rf{vec1-osc} is given by 
\be  \label{osc-exact}
\left( \ba{c}  {\pmb  x}_{n+1} \\  {\pmb p}_{n+1} \ea \right) = \left( 
\ba{cc}   \cos  \Omega  h_n & { \Omega}^{-1} \sin  \Omega  h_n \\ 
-  \Omega \sin  \Omega  h_n  & \cos  \Omega  h_n  \ea \right) 
\left( \ba{c}   {\pmb x}_n \\  {\pmb p}_n \ea \right) + \m 2 \Omega^{-2} \sin^2 \frac{\Omega h_n}{2} \ {\pmb a} \\  \Omega^{-1} \sin \Omega h_n \ {\pmb a}  \ema , 
\ee
where $ h_n$ is an arbitrary variable time step. 
\end{prop}

\begin{cor}[\cite{Ci-oscyl}, Prop.\ 20]   \label{cor-harm}
The formulas \rf{osc-exact} can be rewritten in the following equivalent form:
\be  \ba{l} \label{osc-exact1}
\pmb{\delta}_n^{-1} ({\pmb x}_{n+1} - {\pmb  x}_n ) = \frac{1}{2} \left(   {\pmb p}_{n+1} + {\pmb p}_n \right) \ , \\[2ex]
\pmb{\delta}_n^{-1} ( {\pmb p}_{n+1} - {\pmb p}_n ) = - \frac{1}{2} \Omega^2 ({\pmb x}_{n+1} + {\pmb x}_n ) + {\pmb a} \ .
\ea \ee
where
\be  \label{delta-OM}
   \pmb{\delta}_n = 2 \Omega^{-1} \tan \frac{h_n \Omega}{2} \equiv h_n \tanc \frac{h_n \Omega}{2} 
\ee
\end{cor}

Here and below we use notation \ $\tanc (z) = z^{-1} \tan z$. It is worthwhile to notice that \ $\tanc (z)$ \ is an even function and depends analytically on $z^2$.

As one could expect, the exact discrete harmonic oscillator equations preserve exactly the total energy. 

\begin{prop}[\cite{Ci-oscyl}, Prop.\ 22]  \label{prop-harm}
If \ $\Omega$ is a symmetric matrix ($\Omega^T = \Omega$), then 
\be
I_n : = \frac{1}{2} \scal{{\pmb p}_n}{{\pmb p}_n} + \frac{1}{2} \scal{ {\pmb x}_n}{ \Omega^2  {\pmb x}_n} - \scal{ {\pmb  x}_n }{{\pmb a} } 
\ee 
is an integral of motion (i.e., $I_{n+1} = I_n$) of the discrete multidimensional 
harmonic oscillator equations \rf{osc-exact1}. Here the bracket denotes the scalar product in $\R^m$. 
\end{prop}

Exact discretization seems to be of limited value because, in order to apply it, we need to know the explicit solution of the considered system. However, there exist  non-trivial applications of exact discretizations.  
In the case of the classical Kepler problem we succeeded in using the exact discretization of the harmonic oscillator in two different ways, obtaining two different numerical integrators preserving all trajectories and integrals of motion \cite{Ci-Kep,Ci-Koz,Ci-oscyl}. The exact discretization of the harmonic oscillator equation can also be used to in the integration of some partial differential equations by Fourier transformation \cite{Ci-oscyl}.

\section{Locally exact numerical schemes} 

The central topic of our paper is another fruitful direction of using exact integrators,  namely the so called locally exact discretizations \cite{Ci-oscyl,CR-PRE}. 

\subsection{Motivating example} 

First, we recall our earlier results concerning one-dimensional Hamiltonian systems $\dot p = - V' (x)$, $\dot x = p$, see \cite{CR-long,CR-PRE}. We tested the following class of numerical integrators
\be \ba{l} \label{grad-delV}   \dis
\frac{x_{n+1} - x_n}{\delta_n} = \frac{1}{2} \left( p_{n+1} + p_n \right) \ . \\[3ex] \dis
\frac{p_{n+1} - p_n}{\delta_n} =  - \frac{ V (x_{n+1}) - 
V (x_n) }{x_{n+1} - x_n} \ ,
\ea \ee
where $\delta_n$ is a function defined by
\be \label{deltan} \dis
\delta_n =
\frac{2}{ \omega_n } \tan\frac{ h_n \omega_n  }{2} \ , \qquad 
\omega_n = \sqrt{   V'' (\bar x)  } \ \ ,  
\ee
and, in general,  $\bar x$ may depend on $n$. For simplicity, we formally assume $V'' (\bar x) > 0$. However, in the case of non-positive $V'' (\bar x) $ one can use the same formula,  for details and final results see \cite{CR-PRE}. 
 The simplest choice is $\bar x = x_0$, where $V'(x_0)=0$  (small oscillations around the stable equilibrium). In this case $\delta_n$ does not depend on $n$. The resulting scheme, known as MOD-GR, was first presented in \cite{CR-long}. In  \cite{CR-PRE} we considered the case 
 $\bar x = x_n$ (GR-LEX) and its symmetric (time-reversible) modification \ $\bar x = \frac{1}{2} (x_n + x_{n+1})$ (GR-SLEX). In both cases  $\bar x$ is changed at every step. The names MOD-GR, GR-LEX, GR-SLEX first appeared in \cite{CR-PRE2,CR-BIT}. 

Any numerical scheme from the family defined by \rf{grad-delV} and \rf{deltan} is locally exact (or, more precisely, locally exact at $\bar x$), which means that its linearization at $\bar x$ coincides with the exact discretization of the linearized system. Indeed, the linearized system, given by  $\dot p = - V'(\bar x) - V''(\bar x) \xi$, $\dot \xi = p$, is a particular case of \rf{vec1-osc} and admits the exact discretization \rf{osc-exact1}, \rf{delta-OM}, where we have to identify $\Omega^2 = V'' (\bar x)$, ${\pmb a} = - V' (\bar x)$. This exact discretization coincides with the linearization of \rf{grad-delV}, obtained by substituting $V (x_n) \approx V (\bar x) + V' (\bar x) \xi_n + \frac{1}{2} V'' (\bar x) \xi_n^2$, compare \cite{CR-PRE}.
In this paper we are going to present a generalization of this approach.

\subsection{Local exactness}
\label{sec-locex}

Motivated by the results of \cite{CR-long,CR-PRE} we propose the following definitions.

\begin{Def} \label{def-locex-at}
A numerical scheme ${\pmb x}_{n+1} = \Phi ({\pmb x}_n, h_n)$ for an autonomous equation $\dot {\pmb x} = F ({\pmb x})$  is {\it locally exact} at $\pmb{\bar x}$ if its linearization around $\pmb{\bar x}$ is identical with the exact discretization of the differential equation linearized around $\pmb{\bar x}$. 
\end{Def}

The scheme MOD-GR is locally exact at a stable equilibrium only, GR-LEX is locally exact at $x_n$ (for any $n$), and, finally, GR-SLEX is locally exact at $\frac{1}{2} (x_n + x_{n+1})$ (for any $n$). In the case of implicit numerical schemes the function $\Phi ({\pmb x}_n, h_n)$ is, of course, an implicit function.

\begin{Def}  \label{def-locex}
A numerical scheme ${\pmb x}_{n+1} = \Phi ({\pmb x}_n, h_n)$ for an autonomous equation $\dot {\pmb x} = F ({\pmb x})$  is {\it locally exact} if there exist a sequence $\pmb{\bar x}_n$ such that $\pmb{\bar x}_n  - {\pmb x}_n = O (h_n)$ and the scheme is locally exact at $\pmb{\bar x}_n $ (for any $n$). 
\end{Def}

Therefore GR-LEX and GR-SLEX are locally exact. 
Similar concept (``linearization-preserving'' schemes)  has been recently formulated in \cite{MQTse} (Definition 2.1). An integrator is said to be linearization-preserving if it is linearization preserving at all fixed points. 
This definition is weaker than our Definition~\ref{def-locex}. All locally exact schemes are linearization-preserving but, in general, linearization-preserving scheme has not to be locally exact in our sense. 
For instance, the scheme MOD-GR (see \cite{CR-long,CR-BIT}) is linearization-preserving (provided that $V (x)$ has only one stable equilibrium) but is not locally exact.

\subsection{Exact discretization of linearized equations}
\label{sec-ex-lin}

As an immediate corollary from Definition~\ref{def-locex} we have that the exact discretization of the linearization of a given nonlinear system is locally exact. 
We confine ourselves to autonomous systems of the form
\be  \label{firstord}
  \pmb{\dot x} = F ({\pmb x}) 
\ee
where ${\pmb x} (t) \in \R^d$. 
If ${\pmb x}$ is near a fixed $\pmb{\bar x}$, then \rf{firstord} can be approximated by
\be  \label{approx1}
\dot {\pmb \xi} = F' (\pmb{\bar x}) {\pmb \xi}  +  F (\pmb{\bar x}) 
\ee
where $F'$ is the Fr\'echet derivative (Jacobi matrix) of $F$ and 
\be  \label{xi}
 {\pmb \xi} = {\pmb x} - \pmb{\bar x} \ . 
\ee
The exact discretization of the approximated equation \rf{approx1} is given by
\be  \label{exact-approx}
{\pmb \xi}_{n+1} = e^{ h_n F' (\pmb{\bar x})} {\pmb \xi}_n + \left( e^{ h_n F' (\pmb{\bar x})} - 1 \right) \left( F' (\pmb{\bar x}) \right)^{-1} F (\pmb{\bar x}) \ , 
\ee
provided that \ $\det F' (\pmb{\bar x}) \neq 0$. We make this assumption here and throughout the paper. 
\ods
\begin{cor}  \label{cor-linex}
The exact discretization of the linearization of $\pmb{\dot x} = F (\pmb{x})$ around $\pmb{\bar x}$ is given by: 
\be  \label{lin-ex}
{\pmb x}_{n+1} - \pmb{\bar x} = e^{ h_n F' (\pmb{\bar x})} \left( x_n - \pmb{\bar x} \right) + \left( e^{ h_n F' (\pmb{\bar x})} - 1 \right) \left( F' (\pmb{\bar x}) \right)^{-1} F (\pmb{\bar x})  \ . 
\ee
The scheme \rf{lin-ex} is locally exact. 
\end{cor}
\ods

If we put $\pmb{\bar x} = {\pmb x}_n$ into \rf{lin-ex} (thus changing $\pmb{\bar x}$ at each step), then ${\pmb \xi}_n = 0$ and we obtain:
\be  \label{le-pred}
  {\pmb x}_{n+1} = {\pmb x}_n + \left( e^{  h_n F' ({\pmb x}_n) } - 1 \right) ( F'({\pmb x}_n) )^{-1} F ({\pmb x}_n) 
\ee
This method is well known as the exponential difference equation (see \cite{Pope}, formula (3.4)) or, more recently, as the exponential Euler method \cite{MW}:
\be
{\pmb x}_{n+1} = {\pmb x}_n +  h_n \ \ph_1 (  h_n F'({\pmb x}_n) ) \ F ({\pmb x}_n) \ ,
\ee
where $\ph_1 (z) = z^{-1} (e^z - 1)$ (and $\ph_1 (0) := 1$).  
Another possibility, $\pmb{\bar x} =  \pmb{\bar x}^+_{n}$, where
\be
\pmb{\bar x}^+_{n} = \frac{1}{2} \left( {\pmb x}_{n+1} + \pmb{x}_n \right) \ ,
\ee
leads to 
\be
    {\pmb x}_{n+1} - \pmb{x}_n = 2 ( F' (\pmb{\bar x}^+_{n}))^{-1} \tanh \left( \frac{1}{2}  h_n F'(\pmb{\bar x}^+_{n}) \right) F ( \pmb{\bar x}^+_{n}) \ .
\ee
A third natural possibility, $\pmb{\bar x} = \pmb{x}_{n+1}$, yields the scheme:
\be  \label{imp-Eu1}
  {\pmb x}_{n+1} = {\pmb x}_n + \left( 1 - e^{ - h_n F' ({\pmb x}_{n+1}) } \right) ( F'({\pmb x}_{n+1}) )^{-1} F ({\pmb x}_{n+1}) \ ,
\ee
which may be rewritten in terms of $\ph_1$ as
\be
{\pmb x}_{n+1} = {\pmb x}_n +  h_n \ \ph_1 ( - h_n F'({\pmb x}_{n+1}) ) \ F ({\pmb x}_{n+1}) \ ,
\ee
\ods

\begin{ex}
In the case $\dot p = - V' (x)$, $\dot x = p$, we have 
\be  \label{FF'}
F (\pmb{x}) = \left( \ba{c}   p \\  - V'(x) \ea \right) \ , \qquad F'(\pmb{x}) = \mm 0 & 1 \\ - V'' (x) & 0 \ema  \ ,
\ee
and the scheme \rf{le-pred} assumes the form
\be \ba{l} \dis
x_{n+1} = x_n + \frac{\sin(h_n\omega_n)}{\omega_n}\ p_n - \frac{1-\cos(h_n\omega_n)}{\omega_n^2}\ V'(x_n) \ , \\[3ex]\dis 
p_{n+1} = p_n \cos(h_n\omega_n) -  \frac{\sin(h_n\omega_n)}{\omega_n}\ V'(x_n) \ ,
\ea \ee
where $\omega_n = \sqrt{ V'' (x_n)}$. Global properties of this scheme are rather poor but it can serve as a very good predictor \cite{CR-PRE}. 
\end{ex}

\subsection{Linear stability of locally exact integrators}
\label{sec-stabil} 

Locally exact numerical schemes have excellent qualitative behavior around all fixed points of the considered system. 

\begin{prop}
If the equation $\pmb{\dot x} = F (\pmb{x})$ has a fixed point at \ $\pmb{x} = \pmb{\bar x}$, then 
all its locally exact discretizations  have a fixed point at $\pmb{x}_n = \pmb{\bar x}$, as well.  
\end{prop}

\begin{Proof} If \ $F (\pmb{\bar x}) = 0$, then equation \rf{lin-ex} becomes
\be  \label{linex1}
{\pmb x}_{n+1} - \pmb{\bar x} = e^{ h_n F' (\pmb{\bar x})} \left( x_n - \pmb{\bar x} \right) \ . 
\ee
We require that the scheme $\pmb{x}_{n+1} = \Phi (\pmb{x}_n, h)$ is a locally exact discretization of $F (\pmb{\bar x}) = 0$, i.e., its linearization,  given by
\be
\pmb{x}_{n+1} =  \Phi (\pmb{\bar x}, h) + \Phi' (\pmb{\bar x}, h) (\pmb{x}_n - \pmb{\bar x}) \ ,
\ee
coincides with \rf{linex1}. Hence
\be
 e^{ h_n F' (\pmb{\bar x})} = \Phi' (\pmb{\bar x}, h) \ , \qquad \Phi (\pmb{\bar x}, h) = \pmb{\bar x} \ ,
\ee
which means that $\pmb{\bar x}$ is a fixed point of the system  $\pmb{x}_{n+1} = \Phi (\pmb{x}_n, h)$. 
\end{Proof}

Stability of numerical integrators can be roughly defined as follows: 
``the numerical solution provided by a stable numerical integrator does not tend to infinity when the exact solution is bounded'' (see \cite{CCO}, p. 358). The integrator which is stable when applied to linear equations is said to be linearly stable. We may use the notion of A-stability (see, e.g., \cite{Is-book}): an integrator is said to be A-stable, if discretizations of stable linear equations are stable as well (equation $\dot x = \lambda x$ is stable if $\text{Re} \lambda < 0$). 

Making one more assumption (quite natural, in fact) that the discretization of a linear system is linear, we see that locally exact integrators are linearly stable. Indeed, solutions of any locally exact discretization have the same trajectories as corresponding exact solutions.  

\begin{cor}
Any locally exact numerical scheme is linearly stable and, in particular, A-stable. 
\end{cor}

Some simple examples illustrating this corollary will be given in section~\ref{sec-A-ex}. 
In fact, much stronger result holds: a locally exact discretization yields the best (exact) simulation of a linear equation in the neighborhood of a fixed point.  In particular, locally exact discretizations preserve any qualitative features of trajectories of linear equations (up to round-off errors, of course).  

Numerical experiments show that locally exact schemes are exceptionally stable also for some simple nonlinear systems \cite{CR-long,CR-PRE,CR-BIT}, but we have no theoretical results concerning the stability (e.g., algebraic stability \cite{Is-book}) in the nonlinear case.

\section{Locally exact modifications of popular one-step numerical schemes}
\label{sec-pop}

 We are going to use local exactness as a criterion to select numerical schemes of higher accuracy from a family of non-standard integrators. Our working algorithm to derive such ``locally exact modifications'' of numerical integrators of the form \rf{grad-delH} assumes that $\delta_n$ depends only on \ $\bar x$, $\bar p$ (or, in more general case,  on $\pmb{\bar x}_n$) and $h_n$. We usually consider only three different sequences $\pmb{\bar x}_n $, namely $\pmb{\bar x}_n  = {\pmb x}_n$, $\pmb{\bar x}_n  = {\pmb x}_{n+1}$ and $\pmb{\bar x}_n  = \pmb{x}^+_n \equiv \frac{1}{2} ({\pmb x}_n + {\pmb x}_{n+1})$. The problem of finding the best sequence for a given numerical scheme seems to be interesting but has not been considered yet. 

In this section we apply this  procedure to a number of standard numerical methods. We do not claim that the resulting numerical schemes (sometimes new, sometimes already known) are much better than the corresponding schemes before the modification. Their accuracy should be considerably better but the computing cost is surely much higher. For every scheme separately  one has to check the final outcome and to decide whether the modification is really profitable. Numerical simulations of small oscillations around stable fixed points and the results of \cite{CR-long,CR-PRE} are quite promising. 

\subsection{Locally exact explicit Euler scheme}

The explicit Euler scheme for \rf{firstord} is given by ${\pmb x}_{n+1} = {\pmb x}_n +  h_n F ({\pmb x}_n)$. 
We postulate the following generalization of this scheme:
\be   \label{gen-ex-Euler}
     {\pmb x}_{n+1} = {\pmb x}_n + \pmb{\delta}_n  F ({\pmb x}_n)
\ee
where $\pmb{\delta}_n$ is a variable matrix which is defined by the requirement that 
the linearization of \rf{gen-ex-Euler} coincides with \rf{exact-approx}. We linearize 
\rf{gen-ex-Euler} by substituting \rf{xi} (i.e., ${\pmb x}_n + \pmb{\bar x} + {\pmb \xi}_n$): 
\be  \label{eu-11}
  {\pmb \xi}_{n+1} = {\pmb \xi}_n + \pmb{\delta}_n \left( F (\pmb{\bar x}) + F' (\pmb{\bar x}) {\pmb \xi}_n \right) \ 
\ee
and assume that $\pmb{\delta}_n$ depends only on ${\pmb{\bar x}}$ and $h_n$. Equation~\ref{eu-11} is identical with \rf{exact-approx} iff 
\be  \label{cond-1}
1 + \pmb{\delta}_n F' (\pmb{\bar x}) = e^{ h_n F' (\pmb{\bar x})} , \qquad 
\pmb{\delta}_n = \left( e^{h_n F' (\pmb{\bar x})} - 1 \right) ( F' (\pmb{\bar x}) )^{-1} \ .
\ee
Once can easily see that both equations \rf{cond-1} are equivalent provided that $F' (\pmb{\bar x})$ is non-degenerate.  

Substituting $\pmb{\delta}_n$ from \rf{cond-1} into \rf{gen-ex-Euler} we obtain a class of locally exact modifications of the explicit Euler scheme 
\be  \label{lex-explEuler}
{\pmb x}_{n+1} = {\pmb x}_n + \big( e^{ h_n F' (\pmb{\bar x})} - 1 \big) ( F' (\pmb{\bar x}) )^{-1} F ({\pmb x}_n) \ .
\ee
It seems that the choice $\pmb{\bar x} = {\pmb x}_n$ is most natural in that case (then, by the way, \rf{lex-explEuler} coincides with \rf{le-pred}). 

\begin{cor}
Locally exact modification of the explicit Euler scheme is given by
\be  \label{lee-Euler}
{\pmb x}_{n+1} = {\pmb x}_n + \big( e^{ h_n F' ({\pmb x}_n)} - 1 \big) ( F' ({\pmb x}_n) )^{-1} F ({\pmb x}_n)
\ee
\end{cor}

\subsection{Locally exact implicit Euler scheme}

The implicit Euler scheme for \rf{firstord} is given by ${\pmb x}_{n+1} = {\pmb x}_n +  h_n F ({\pmb x}_{n+1})$. 
We postulate the following generalization of this scheme:
\be   \label{gen-imp-Euler}
     {\pmb x}_{n+1} = {\pmb x}_n + \pmb{\delta}_n  F ({\pmb x}_{n+1})
\ee
where $\pmb{\delta}_n$ is a variable matrix which is defined by the requirement that 
the linearization of \rf{gen-imp-Euler} coincides with \rf{exact-approx}. We linearize 
\rf{gen-imp-Euler} by substituting \rf{xi} (i.e., ${\pmb x}_n + \pmb{\bar x} + {\pmb \xi}_n$): 
\be
  {\pmb \xi}_{n+1} = {\pmb \xi}_n + \pmb{\delta}_n \left( F (\pmb{\bar x}) + F' (\pmb{\bar x}) {\pmb \xi}_{n+1}\right) \ ,
\ee
i.e., $\big( 1 -\pmb{\delta}_n F'(\pmb{\bar x}) \big) {\pmb \xi}_{n+1} = {\pmb \xi}_n + \pmb{\delta}_n F (\pmb{\bar x})$, which is identical with \rf{exact-approx} iff 
\be  \label{cond-2} \ba{l} 
\big( 1 - \pmb{\delta}_n F' (\pmb{\bar x}) \big)^{-1}  = e^{ h_n F' (\pmb{\bar x})} \\[2ex] 
\big( 1 - \pmb{\delta}_n F' (\pmb{\bar x}) \big)^{-1} \pmb{\delta}_n = \big( e^{ h_n F' (\pmb{\bar x})} - 1 \big) ( F' (\pmb{\bar x}) )^{-1} 
.
\ea \ee
We easily see that $\pmb{\delta}_n$ can be computed independently  from any of these two equations. Fortunately, the result is the same: 
\be
 \pmb{\delta}_n = \big(  1 - e^{-  h_n F' (\pmb{\bar x})} \big) ( F' (\pmb{\bar x}) )^{-1} 
\ee
Substituting $\pmb{\delta}_n$ from \rf{cond-2} into \rf{gen-imp-Euler} we get
\be  \label{lex-impEuler-gen}
{\pmb x}_{n+1} = {\pmb x}_n + \big(  1 - e^{-  h_n F' (\pmb{\bar x})} \big) ( F' (\pmb{\bar x}) )^{-1} F ({\pmb x}_{n+1}) \ ,
\ee
or, in an equivalent way,
\be
{\pmb x}_{n+1} = {\pmb x}_n +  h_n \ \ph_1 ( -  h_n F'(\pmb{\bar x}) ) \ F ({\pmb x}_{n+1}) \ .
\ee
In this case it is not clear what is the most natural identification. We can choose either $\pmb{\bar x} = {\pmb x}_n$, or $\pmb{\bar x} = {\pmb x}_{n+1}$.  

\begin{cor}
Locally exact modification of the implicit Euler scheme is given either by
\be  \label{leie-Euler}
{\pmb x}_{n+1} = {\pmb x}_n + \big(  1 - e^{-  h_n F' ({\pmb x}_n)} \big) ( F' ({\pmb x}_n) )^{-1} F ({\pmb x}_{n+1}) \ ,
\ee
or, by
\be
\label{leii-Euler}
{\pmb x}_{n+1} = {\pmb x}_n + \big(  1 - e^{-  h_n F' ({\pmb x}_{n+1})} \big) ( F' ({\pmb x}_{n+1}) )^{-1} F ({\pmb x}_{n+1}) \ . 
\ee
\end{cor}

\subsection{Locally exact implicit midpoint rule} 

We postulate the following extension of the implicit midpoint rule: 
\be  \label{imr-del}
{\pmb x}_{n+1} - {\pmb x}_n  = \pmb{\delta}_n  F \left(\frac{{\pmb x}_{n+1} + {\pmb x}_n}{2} \right) \ ,
\ee
where $\pmb{\delta}_n$ depends, as usual, on $\pmb{\bar x}$ and $h_n$. The scheme \rf{imr-del} is locally exact for
\be  \label{imr-deltan}
\pmb{\delta}_n  = 2 ( F' (\pmb{\bar x}) )^{-1} \tanh \left( \frac{h_n F'(\pmb{\bar x})}{2} \right)   \ , 
\ee
which will be shown, in more general framework, in section~\ref{le-gen}.  
In some applications the Taylor series may be useful:
\be
\pmb{\delta}_n = 1 - \frac{1}{12}  h_n^2 (F'(\pmb{\bar x}))^2 + \frac{1}{120}  h_n^4 (F'(\pmb{\bar x}))^4 + \ldots  
\ee
Thus locally exact modification of the implicit midpoint rule reads
\be  \label{le-midpoint}
{\pmb x}_{n+1} - {\pmb x}_n  = h_n \left( \tanhc \frac{ h_n F'(\pmb{\bar x})}{2} \right) \ F \left( \frac{{\pmb x}_{n+1} + {\pmb x}_n}{2} \right) 
\ee
The most natural choice of $\pmb{\bar x}$ seems to be at the midpoint, $\pmb{\bar x} = \frac{1}{2} \left( {\pmb x}_n + {\pmb x}_{n+1} \right)$. However, in order to diminish the computation cost, the choice $\pmb{\bar x} = {\pmb x}_n$ can also be considered (because then the Jacobian $F' (\pmb{\bar x})$ is evaluated outside iteration loops).

\subsection{Locally exact trapezoidal rule} 

We postulate the following extension of the trapezoidal  rule: 
\be  \label{tr-del}
 {\pmb x}_{n+1} - {\pmb x}_n = \pmb{\delta}_n \ \frac{ F ({\pmb x}_{n+1}) + F ({\pmb x}_n) }{2} \ ,
\ee
Linearizing \rf{tr-del} around $\pmb{\bar x}$ we get 
\be
  {\pmb \xi}_{n+1} - {\pmb \xi}_n = \pmb{\delta}_n \left( F (\pmb{\bar x}) + F'(\pmb{\bar x}) \ \frac{ {\pmb \xi}_{n+1} + {\pmb \xi}_n}{2} \right) 
\ee
In section~\ref{le-gen} we will show that the scheme \rf{tr-del} is locally exact for
\be
\pmb{\delta}_n  = 2 ( F' (\pmb{\bar x}) )^{-1} \tanh \left( \frac{h_n F'(\pmb{\bar x})}{2} \right) \ .  
\ee
Thus locally exact modification of the trapezoidal rule 
is given by
\be  \label{le-trapez}
{\pmb x}_{n+1} - {\pmb x}_n  =  h_n \left( \tanhc \frac{ h_n F'(\pmb{\bar x})}{2} \right) \ \frac{ F ({\pmb x}_{n+1}) + F ({\pmb x}_n) }{2} \ . 
\ee
There are two natural choices of $\pmb{\bar x}$. Either (in order  to minimize the computational costs) we can take $\pmb{\bar x} = {\pmb x}_n$, or (in order to obtain a time-reversible scheme) we can take $\pmb{\bar x} = \frac{1}{2} \left( {\pmb x}_n + {\pmb x}_{n+1} \right)$.

\subsection{A-stability of locally exact modifications}
\label{sec-A-ex} 

In order to illustrate general results of section~\ref{sec-stabil}, we 
will apply four schemes presented above to the one-dimensional linear equation $\dot x = \lambda x$. 
The locally exact explicit Euler scheme 
yields
\be
  x_{n+1} - x_n = \left( e^{h_n \lambda} - 1 \right)  x_n \ . 
\ee
The locally exact implicit Euler scheme yields
\be
  x_{n+1} - x_n = \left( 1 - e^{- h_n \lambda} - 1 \right)  x_{n+1} \ , 
\ee
Both resulting equations are identical: $x_{n+1} = e^{h_n \lambda} x_n$ and yield the exact discretization.   
The implicit midpoint and trapezoidal rules yield an identical equation, namely 
\be  \label{ideq} 
 x_{n+1} - x_n = \left( \tanh \frac{h_n \lambda}{2} \right) \left( x_{n+1} + x_n \right) \ .
\ee
Computing $x_{n+1}$ from \rf{ideq}  we get the exact discretization  $x_{n+1} = e^{h_n \lambda} x_n$ again. 
In particular, if $\text{Re} \lambda < 0$, then $x_n \rightarrow 0$ for $n \rightarrow \infty$ (for any constant time step $h_n = h = \const$). A-stability is evident in all these cases.

\subsection{A large class of locally exact integrators}
\label{le-gen}

All numerical schemes presented above are particular cases of  the following class of numerical schemes for the equation ${\pmb{\dot x}} = F ({\pmb x})$: 
\be \label{scheme-gen}
  {\pmb x}_{n+1} - {\pmb x}_n = \pmb{\delta} ({\pmb{\bar x}}) \Psi ({\pmb x}_n, {\pmb x}_{n+1} ) \ .
\ee
We assume the consistency conditions:
\be  \label{consis}
  \Psi ({\pmb x}, {\pmb x}) = F ({\pmb x}) \ , \qquad  F'({\pmb x}) = \Psi_1 ({\pmb x}, {\pmb x}) + \Psi_2 ({\pmb x}, {\pmb x})  \ , 
\ee
where $\Psi_1, \Psi_2$ are partial Fr\'echet derivative with respect to the first and second vector variable, respectively (thus $\Psi_1, \Psi_2$ are $d\times d$ matrices). We also denote
\be  \label{xx}
 {\bar \Psi} = \Psi ({\pmb {\bar x}}, {\pmb {\bar x}}) \ , \quad {\bar \Psi}_1 =  \Psi_1 ({ \pmb {\bar x}}, {\pmb {\bar x}}) \ , \quad {\bar \Psi}_2 =  \Psi_2 ({ \pmb {\bar x}}, {\pmb {\bar x}}) \ .
\ee

\begin{prop}  \label{prop-legen}
The numerical scheme \rf{scheme-gen}, where $\Psi$ satisfies \rf{consis}, is locally exact for
\be  \label{delta-gen}
  \pmb{\delta} ({\pmb{\bar x}}) = \left( e^{h_n F' (\pmb{\bar x})} - 1 \right) \left( F'(\pmb{\bar x}) + {\bar \Psi}_2 \left( e^{h_n F' (\pmb{\bar x})} - 1 \right)  \right)^{-1} \ ,
\ee
\end{prop}

\begin{Proof} The exact discretization of the linearization of equation ${\pmb{\dot x}} = F ({\pmb x})$ is given by  \rf{exact-approx}. The linearization of the scheme \rf{scheme-gen} (at $\pmb{x}_n = \pmb{\bar x}$)  reads
\be  \label{ksi1}
  {\pmb \xi}_{n+1} - {\pmb  \xi}_n = \pmb{\delta} ( {\bar \Psi}_1 {\pmb \xi}_n + {\bar \Psi}_2 {\pmb \xi}_{n+1} ) + \pmb{\delta} {\bar \Psi}
\ee
where  $\pmb{\delta} = \pmb{\delta} (\pmb{\bar x})$ and we use \rf{xx}. Identifying \rf{ksi1} with \rf{exact-approx} we get a system of two equations: 
\be \label{fir}
(1 - \pmb{\delta} {\bar \Psi}_2)^{-1} (1 + \pmb{\delta} {\bar \Psi}_1) = e^{h_n F'} \ , 
\ee
\be  \label{sec}
(1 - \pmb{\delta} {\bar \Psi}_2)^{-1} \pmb{\delta} {\bar \Psi} = \left( e^{h_n F'} - 1 \right) (F')^{-1} F \ ,
\ee
where $F = F (\pmb {\bar x})$ and $F'= F' (\pmb {\bar x})$. 
Equation \rf{fir} implies
\be  \label{symm}
\pmb{\delta} \left( {\bar \Psi}_1 + {\bar \Psi}_2  e^{h_n F'} \right) = e^{h_n F'} - 1 
\ee
Eliminating ${\bar \Psi}_1$ from \rf{symm} (by virtue of \rf{consis}), we get
\be
\pmb{\delta} {\bar \Psi}_2 \left( e^{h_n F'} - 1 \right) + \pmb{\delta} F' = e^{h_n F'} - 1 
\ee
which implies \rf{sec} (note that $F = {\bar \Psi}$). Therefore, the local exactness imposes  only one condition for $\pmb{\delta}$, namely \rf{symm}. 
\end{Proof}

Expressing $\pmb{\delta}$ given by \rf{delta-gen} in a more symmetric way, we obtain another, equivalent, form of locally exact scheme \rf{scheme-gen}:  
\be
 {\pmb x}_{n+1} - {\pmb x}_n =  h_n  \tanhc \frac{h_n F'}{2}  \left( 1 +  h_n \left( {\bar \Psi}_2 - {\bar \Psi}_1 \right) \tanhc \frac{h_n F'}{2}  \right)^{-1}   \Psi ({\pmb x}_n, {\pmb x}_{n+1} ) .
\ee

All numerical schemes presented above can be considered as particular cases of \rf{scheme-gen}, namely: \par
\renewcommand{\arraystretch}{1.5}\par
\begin{tabular}{lll}
$\bullet$ &  \text{explicit Euler scheme:} & $\Psi ({\pmb x}_n, {\pmb x}_{n+1} ) = F (\pmb{x}_n)$ \ , \\ 
$\bullet$ & \text{implicit Euler scheme:} & $\Psi ({\pmb x}_n, {\pmb x}_{n+1} ) = F (\pmb{x}_{n+1})$ \ , \\ 
$\bullet$ & \text{implicit midpoint rule:} & $\Psi ({\pmb x}_n, {\pmb x}_{n+1} ) = F (\pmb{x}^+_n)$ \ , \\
$\bullet$ & \text{trapezoidal rule:} & $\Psi ({\pmb x}_n, {\pmb x}_{n+1} ) = \frac{1}{2} \left( F (\pmb{x}_n) + F(\pmb{x}_{n+1}) \right)$ \ .
\end{tabular} \par
\renewcommand{\arraystretch}{1} \ods
\no Therefore, the proof of Proposition~\ref{prop-legen} is valid for all these cases.

\section{Locally exact discrete gradient schemes for one-dimensional Hamiltonian systems}
\label{sec-locex-one}

The discrete gradient scheme (or modified midpoint rule) is a conservative integrator which has been used since many years for simulating dynamics of systems of particles  \cite{Gree,LaG}. 
More recently discrete gradient methods have been developed in the context of geometric numerical integration \cite{MQ}, see \cite{Gon,IA}.   In particular, Quispel and his co-workers constructed numerical integrators preserving integrals of motion of a given system of ordinary differential equations \cite{MQR2,QC,QT}. 
In this section we consider  Hamiltonian systems with one degree of freedom: 
\be  \label{Ham} 
  \dot x = H_p \ , \qquad \dot p = - H_x \ , 
\ee
where $H = H (x, p)$ is a given function (sufficiently smooth), subscripts denote partial differentiation and the dot denotes the total derivative respect to $t$. The Hamiltonian $H (x, p)$ is an integral of motion (the energy integral). 
In order to make this paper more self-contained, in section~\ref{sec-grad1} we present a locally exact symmetric discrete gradient scheme, first obtained in \cite{CR-BIT}. In section~\ref{sec-one-incre} we derive new locally exact scheme modifying the coordinate increment discrete gradient method.

\subsection{Locally exact symmetric discrete gradient scheme}
\label{sec-grad1}

Following \cite{CR-BIT}, 
we consider a class of non-standard (compare \cite{Mic}) discrete gradient schemes for the system \rf{Ham}:  
\be  \ba{l} \dis \label{grad-delH} 
 \frac{x_{n+1} - x_n}{\delta_n} =  \frac{  H (x_{n+1}, p_{n+1}) + H (x_n, p_{n+1}) - H (x_{n+1}, p_n) -  H (x_n, p_n) }{2 (p_{n+1}- p_n) } \ , 
\\[4ex] \dis
\frac{p_{n+1} - p_n}{\delta_n} =  \frac{  H (x_n, p_{n+1}) + H (x_n, p_n) - H (x_{n+1}, p_{n+1}) -  H (x_{n+1}, p_n) }{2 (x_{n+1}- x_n) } \ ,   
\ea \ee
where $\delta_n$ is an arbitrary positive function of $h_n, x_n, p_n, x_{n+1}, p_{n+1}$ etc. (the time step is denoted by $h_n$). The subscript $n$ indicates that $\delta_n$ may depend on the step $n$.  The discrete system \rf{grad-delH} is manifestly symmetric (time-reversible).

\begin{lem}
The scheme \rf{grad-delH} exactly preserves the energy integral for any $\delta_n$, i.e., $H (x_{n+1}, p_{n+1}) = H (x_n, p_n)$. 
\end{lem}

\begin{Proof} In order to obtain the energy conservation law it is enough to multiply the first equation by $2 (p_{n+1} - p_n)$ and the second equation by $2(x_{n+1} - x_n)$, and to subtract resulting equations. 
\end{Proof}

\begin{prop}   \label{prop-locex}
The discrete gradient scheme \rf{grad-delH} with 
\be  \label{delta-omH} 
  \delta_n = \frac{2}{\omega_n} \tan \frac{h_n \omega_n}{2} \ , \qquad \omega_n = \sqrt{ H_{xx} H_{pp} - H_{xp}^2 } \ , 
\ee
(where $\omega_n$ is evaluated at $\bar x, \bar p$) is locally exact. 
\end{prop}

\begin{Proof} We linearize \rf{Ham}, substituting $x = {\bar x} + \xi$, $p = {\bar p} + \eta$: 
\be  \label{linear-H}
   \dot \xi = H_p + H_{px}  \xi  + H_{pp} \eta \ , \qquad  {\dot \eta} = - H_x -  H_{xx} \xi - H_{xp} \eta \ . 
\ee
The exact discretization of the system \rf{linear-H} is given by 
\be  \label{dis-Ab} 
  \left( \ba{c} \xi_{n+1} \\ \eta_{n+1} \ea \right)  = e^{h_n F'} \left( \ba{c} \xi_n \\ \eta_n \ea \right)   +  \left( e^{h_n F'} - 1 \right) (F')^{-1} F \ ,
\ee 
(compare Proposition~\ref{prop-exact}), where
\be  \label{def-Ab} 
 F = \m H_p \\ - H_x \ema \ , \qquad  F' =  \left( \ba{cc}  H_{xp} & H_{pp} \\ - H_{xx} & - H_{xp} \ea \right) \ .
\ee
Then, we linearize the system \rf{grad-delH} around $\bar x, \bar p$, obtaining
\be  \ba{l} \dis  \label{hlin}
 \frac{\xi_{n+1} - \xi_n}{\delta_n} = H_p + \frac{1}{2} H_{xp} \left( \xi_n + \xi_{n+1} \right) + \frac{1}{2} H_{pp} \left( \eta_n + \eta_{n+1} \right) \ , \\[3ex]  \dis
 \frac{\eta_{n+1} - \eta_n}{\delta_n} = - H_x - \frac{1}{2} H_{xx} \left( \xi_n + \xi_{n+1} \right) - \frac{1}{2} H_{xp} \left( \eta_n + \eta_{n+1} \right) \ ,
\ea \ee
where $x_n = {\bar x} + \xi_n$, \ $p_n = {\bar p} + \eta_n$ and 
 partial derivatives $H_x, H_p, H_{xx}, H_{xp}$ and $H_{pp}$ are evaluated at ${\bar x}, {\bar p}$. 
The system \rf{hlin} can be rewritten in the matrix form
\be  \label{lin-grad}
 \left( 1 - \frac{1}{2} \delta_n F' \right) \m \xi_{n+1} \\ \eta_{n+1} \ema  = \left(1 + \frac{1}{2} \delta_n F' \right) \m  \xi_n \\ \eta_n \ema + \delta_n F  \ ,
\ee
where $F$ and $F'$ are defined by \rf{def-Ab}. 
We easily verify that 
\be \label{A2}
  (F')^2 = - \omega_n^2  \ , \qquad \omega_n^2 = H_{xx} H_{pp} - H_{xp}^2 \ .
\ee
where here (and in many other places) we omit the unit matrix (i.e., we write $\omega_n^2$ instead of $\omega_n^2 I$, etc.). 
Comparing \rf{lin-grad} and \rf{dis-Ab} we obtain   local exactness conditions:  
\be  \ba{l} 
e^{h_n F'} = (1 - \frac{1}{2} \delta_n F')^{-1}(1 + \frac{1}{2} \delta_n F')  \ , \\[3ex] 
( e^{h_n F'} - 1) (F')^{-1} F = (1 - \frac{1}{2} \delta_n F' )^{-1} \delta_n F \ .
\ea \ee 
Substituting the first equation into the second one, we get an identity. Therefore, it is enough to consider the first equation:
\[
  e^{h_n F'} - \frac{1}{2} \delta_n F' e^{h_n F'} = 
1 + \frac{1}{2} \delta_n F'   \quad \Longrightarrow \quad 
\delta_n = 2 (F')^{-1} \tanh \frac{1}{2} h_n F' \ . 
\]
Therefore, $\delta_n$ depends analytically on $(F')^2$ and, by virtue of \rf{A2}, $\delta_n$ is proportional to the unit matrix. Hence $\delta_n$ is indeed a scalar function, given by  
\[
   \delta_n = 2 (F')^{-1} \tanh \frac{1}{2} h_n F' = \frac{2}{\omega_n} \tan \frac{1}{2} h_n \omega_n \ ,
\]
which ends the proof. 
\end{Proof}

As usual, we may take either $\bar x = \bar x_n^+$, $\bar p = \bar p_n^+$ (to keep the scheme symmetric), or $\bar x = x_n$, $\bar p = p_n$ (to minimize the computational cost).  

We point out that the formula \rf{delta-omH} implies some limitations on $h_n$ in the case $\omega_n \in \R$. Certainly we have to require $h_n  \omega_n \neq  \pi + 2 \pi M$ ($M \in \N$), or even $h_n \omega_n <  \pi$.  The last inequality is quite reasonable because it means that $h_n < \frac{1}{2} T_n$, where $T_n = 2\pi/\omega_n$ is a corresponding period, compare \cite{CR-PRE}.

\subsection{Locally exact coordinate increment discrete gradient scheme} 
\label{sec-one-incre}

In the previous section we used symmetric form of the discrete gradient. The case of coordinate increment gradient (see \cite{IA}), although of simpler form,  turned out to be more difficult in the context of locally exact modifications. However, we succeeded to derive such modification also in this case. 

We consider the following non-standard numerical integrator for the Ha\-miltonian system \rf{Ham}: 
\be \ba{l}\dis  \label{incre}
 \frac{x_{n+1} - x_n}{\delta_n} = \frac{H (x_{n+1}, p_{n+1}) - H ( x_{n+1}, p_n)}{p_{n+1} - p_n} \ , \\[2ex]\dis
\frac{p_{n+1} - p_n}{\delta_n} = \frac{ H (x_n, p_n) - H (x_{n+1}, p_n)}{x_{n+1} - x_n} \ , 
\ea \ee
where, similarly as in the formula \rf{grad-delH}, $\delta_n$ is an arbitrary function. 
\ods

\begin{lem}
The scheme \rf{incre} exactly preserves the energy integral for any $\delta_n$, i.e., $H (x_{n+1}, p_{n+1}) = H (x_n, p_n)$. 
\end{lem}
\begin{Proof} We multiply the first equation by $p_{n+1} - p_n$ and the second equation by $x_{n+1} - x_n$. Then, we subtract resulting equations. 
\end{Proof}

\begin{prop}
The scheme \rf{incre} is locally exact for $\delta_n$ given by
\be \label{delta-incre}
\delta_n = \frac{2}{\omega_n \cot \frac{\omega_n h_n}{2} + H_{x p} }  \ , \qquad \omega_n = \sqrt{ H_{xx} H_{pp} - H_{xp}^2 } \ , 
\ee
where derivatives of $H$ are evaluated at $\bar x, \bar p$. 
\end{prop}

\begin{Proof}
Using notation from section~\ref{sec-grad1} we linearize \rf{incre} around $(\bar x, \bar p)$:
\be \ba{l}
\xi_{n+1} - \xi_n = \delta_n H_p + \frac{1}{2} \delta_n H_{pp} (\eta_n + \eta_{n+1}) + H_{xp} \xi_{n+1} \ , \\[2ex]
\eta_{n+1} - \eta_n = - \delta_n H_x - \frac{1}{2} \delta_n H_{xx} (\xi_n + \xi_{n+1}) - H_{xp} \eta_n \ , 
\ea \ee  
which can be rewritten as (compare \rf{def-Ab}): 
\be  \label{xx-1}
\left( 1 - \frac{1}{2} H_{xp} \delta_n - \frac{1}{2} \delta_n F' \right) \m \xi_{n+1} \\ \eta_{n+1} \ema = \left(  1 - \frac{1}{2} H_{xp} \delta_n + \frac{1}{2} \delta_n F' \right) \m \xi_n \\ \eta_n \ema + \delta_n F  .
\ee
Requiring that \rf{xx-1} is identical with the exact discretization \rf{dis-Ab} we get:
\be  \ba{l} \dis  \label{le-cond}
\left( 1 - \frac{1}{2} H_{xp} \delta_n - \frac{1}{2} \delta_n F' \right) e^{h_n F'} = 1 - \frac{1}{2} H_{xp} \delta_n + \frac{1}{2} \delta_n F' \ , \\[2ex]\dis
\left( 1 - \frac{1}{2} H_{xp} \delta_n - \frac{1}{2} \delta_n F' \right)  \left( e^{h_n F'} - 1 \right) (F')^{-1} F = \delta_n F \ .
\ea \ee
Substituting the left-hand side of the first equation into the second equation, we get an identity. We use the first equation to compute $\delta_n$: 
\be
\delta_n = 2 \left( e^{h_n F'} - 1 \right) \left(  ( e^{h_n F'} + 1 ) F' + ( e^{h_n F'} - 1 ) H_{xp} \right)^{-1} \ , 
\ee
which reduces to 
\be
 \delta_n = 2 \left( H_{xp} + F' \cot \frac{h_n F'}{2} \right)^{-1} \ .
\ee
Finally, we observe that function $f (z) = z \cot z$ is  even and depends analytically on $z^2$. Moreover, $(F')^2= -\omega_n^2$ is proportional to the unit matrix. Therefore $\delta_n$ is a scalar function (which is compatible with our assumption), given by \rf{delta-incre}. 
\end{Proof}

\subsection{One-dimensional separable Hamiltonian systems}
\label{sec-one-sep}

 As a simple example we consider the case defined by $H_{xp} = 0$, i.e, 
\be \label{ham-sep}
  H = T (p) + V (x) \ .
\ee
In this case the coordinate increment discrete gradient becomes identical with the symmetric discrete gradient and   
 the scheme \rf{grad-delH} reduces to   
\be  \ba{l} \dis  \label{grad-delTV} 
\frac{x_{n+1} - x_n}{\delta_n} =  \frac{  T  (p_{n+1}) - T (p_n) }{p_{n+1}- p_n } \ ,  
\\[4ex] \dis
\frac{p_{n+1} - p_n}{\delta_n} =  - \frac{  V (x_{n+1}) - V (x_n)}{ x_{n+1}- x_n } \ .
\ea \ee
Local exactness yields 
\be
  \delta_n = h_n \tanc \frac{h_n \omega_n}{2} \ , \qquad \omega_n = \sqrt{ T'' (\bar p) V'' (\bar q) } \ . 
\ee
Separable Hamiltonians  \rf{ham-sep} are 
associated with  many  mechanical systems with one degree of freedom. The case $T(p) = \frac{1}{2} p^2$ was  considered in previous papers \cite{CR-long,CR-PRE,CR-BIT}, where many numerical experiments were reported. 
Assuming $\bar x = x_n$, $\bar p = p_n$  we get a scheme called GR-LEX, while $\bar x = \frac{1}{2} \left( x_n + x_{n+1} \right)$, $\bar p = \frac{1}{2} \left( p_n + p_{n+1} \right)$ yields GR-SLEX \cite{CR-BIT}. The system \rf{Ham} is symmetric (time-reversible).  GR-SLEX preserves this property, while GR-LEX does not.

The discrete gradient schemes GR and MOD-GR are of second order. Locally exact discrete gradient schemes have higher order: GR-LEX is of 3rd order and GR-SLEX is of 4th order, see \cite{CR-PRE}. Numerical experiments presented in \cite{CR-PRE,CR-BIT} have shown that the accuracy of GR-LEX and GR-SLEX is higher by several orders of magnitude when compared with the standard discrete gradient method (while the computational cost is higher at most several times, usually much less). Locally exact modifications turn out to be of considerable advantage.

\section{Locally exact discrete gradient schemes for multidimensional Hamiltonian systems}
\label{sec-locex-multi}

This section contains main results. We extend results of section~\ref{sec-locex-one} constructing  
energy-preserving locally exact discrete gradient schemes for arbitrary multidimensional Hamiltonian systems in canonical coordinates:
\be \label{multiham}
  {\dot x}^k = \frac{\partial H}{\partial p^k} \ , \quad {\dot p}^k = - \frac{\partial H}{\partial x^k} \ .
\ee
where $k = 1,\ldots,m$. 
We obtain two different numerical schemes using either  symmetric discrete gradient or coordinate increment discrete gradient. 

\subsection{Linearization of Hamiltonian systems}

We denote $\pmb{x} = (x^1,\ldots, x^m)^T$, $\pmb{p} = (p^1,\ldots, p^m)$, etc. 
The linearization of \rf{multiham} around \ $\pmb{\bar x}, \pmb{\bar p}$ \ is given by:
\be \ba{l} \dis 
{\dot \xi}^i = H_{p^i} + \sum_{k=1}^m H_{p^i x^k} \xi^k + \sum_{k=1}^m H_{p^i p^k} \eta^k \ , \\[3ex] \dis
{\dot \eta}^i = - H_{x^i} - \sum_{k=1}^m H_{x^i p^k} \eta^k - \sum_{k=1}^m H_{x^i x^k} \xi^k \ ,
\ea \ee
or, in a matrix notation
\be \ba{l}  \label{lin-mham}
 \pmb{\dot \xi} = H_{\pmb{p}} + H_{\pmb{p} \pmb{x}} \pmb{\xi} + H_{\pmb{p} \pmb{p}} \pmb{\eta} \ , \\[2ex]
\pmb{\dot \eta} = - H_{\pmb{x}} - H_{\pmb{x} \pmb{x}} \pmb{\xi} - H_{\pmb{x} \pmb{p}} \pmb{\eta} \ , 
\ea \ee
where derivatives of $H$ are evaluated at $\pmb{\bar x}, \pmb{\bar p}$.  
Note that $m\times m$ matrices $H_{\pmb{x} \pmb{x}}$, $H_{\pmb{x} \pmb{p}}$, $H_{\pmb{p} \pmb{x}}$, $H_{\pmb{p} \pmb{p}}$ satisfy 
\be
 H_{\pmb{p} \pmb{p}}^T = H_{\pmb{p} \pmb{p}} \ , \quad H_{\pmb{x} \pmb{x}}^T = H_{\pmb{x} \pmb{x}}  \ , \quad H_{\pmb{x} \pmb{p}}^T = H_{\pmb{p} \pmb{x}}  
\ee
Equations \rf{lin-mham} can be rewritten also as
\be \label{lin-mmh}
\frac{d}{d t} \m \pmb{\xi} \\ \pmb{\eta} \ema = F' \m \pmb{\xi} \\ \pmb{\eta} \ema  + F 
\ee
where 
\be  \label{F-F'}
 F = \m H_{\pmb{p}} \\ - H_{\pmb{x}} \ema \ , \quad F' = \mm H_{\pmb{p}\pmb{x}} & H_{\pmb{p}\pmb{p}} \\ - H_{\pmb{x}\pmb{x}} & - H_{\pmb{x}\pmb{p}} \ema \ .
\ee

\begin{cor}  \label{cor-linham} 
The exact discretization of linearized Hamiltonian equations \rf{lin-mmh} is given by 
\be  \label{lin-hex}
 \m \pmb{\xi}_{n+1} \\ \pmb{\eta}_{n+1} \ema = e^{h_n F'} \m \pmb{\xi}_n \\ \pmb{\eta}_n \ema + \left( e^{h_n F'} - 1 \right) (F')^{-1} F \ .
\ee
\end{cor}

\begin{Proof}
 The exact discretization of \rf{lin-mmh} is given by \rf{lin-hex} which follows immediately from Corollary~\ref{cor-linex}. 
\end{Proof}

\subsection{Discrete gradients in the multidimensional case} 

Considering multidimensional Hamiltonian systems we will denote 
\be  \label{y=xp}
\pmb{y} := \m \pmb{x} \\  \pmb{p} \ema \ , \quad   \pmb{y}_n := \m \pmb{x}_n \\ \pmb{p}_n \ema \ ,  
\ee
where $\pmb{x} \in \R^m$, $\pmb{p} \in \R^m$, $\pmb{y} \in \R^{2 m}$, etc. In other words, 
\[ 
y^1 = x^1$, $y^2 = x^2$, \ldots,$y^m = x^m$, $y^{m+1} = p^1$,\ldots, $y^{2m} = p^m \ . 
\]
A discrete gradient, denoted by $\bar \nabla H (\pmb{y}_n, \pmb{y}_{n+1})$,  
\be
 {\bar \nabla} H   = \left( \frac{ {\Delta} H}{\Delta y^1}, \frac{ {\Delta} H}{\Delta y^2}, \ldots,\frac{ {\Delta} H}{\Delta y^{2 m}}  \right) \equiv \left( \frac{ \Delta H}{\Delta \pmb{x}}, \frac{ \Delta H }{\Delta \pmb{p}} \right) 
\ee
is defined as an $\R^{2 m}$-valued  function of $\pmb{y}_n,  \pmb{y}_{n+1}$ such that \cite{Gon,MQR2}: 
\be \ba{l}\dis \label{disgrad-def}
 \sum_{k=1}^{2 m}   \frac{ {\Delta} H }{\Delta y^k } \left( y_{n+1}^k - y_n^k  \right)  = H (\pmb{y}_{n+1}) - H (\pmb{y}_n) \ , 
\\[4ex]\dis
{\bar \nabla} H (\pmb{y}, \pmb{y})  = \left( H_{\pmb{x}}, H_{\pmb{p}} \right) \ , 
\ea \ee
where $H_{\pmb{x}}, H_{\pmb{p}}$ are evaluated at $\pmb{y} = (\pmb{x}, \pmb{p})$ and we define
\be  \label{grad-lim}
{\bar \nabla} H (\pmb{y}, \pmb{y}) = \lim_{\pmb{\tilde y} \rightarrow \pmb{y}} {\bar \nabla} H (\pmb{y}, \pmb{\tilde y})
\ee 
when necessary. 

Discrete gradients are non-unique, compare \cite{Gon,IA,MQR2}. Here we confine ourselves to the simplest form of the discrete gradient, namely coordinate increment discrete gradient \cite{IA} and to its symmetrization. 
The coordinate increment discrete gradient is defined by:
\be \ba{l} \dis \label{grad-incre} 
 \frac{\Delta H}{\Delta y^1} = \frac{H(y_{n+1}^1, y_n^2, y_n^3,\ldots, y_n^{2 m}) - H (y_n^1, y_n^2, y_n^3,\ldots, y_n^{2 m})}{y_{n+1}^1 - y_n^1}  , \\[3ex]\dis  
 \frac{\Delta H}{\Delta y^2} = \frac{H(y_{n+1}^1, y_{n+1}^2, y_n^3,\ldots, y_n^{2 m}) - H (y_{n+1}^1, y_n^2, \ldots, y_n^{2 m})}{y_{n+1}^2 - y_n^2}  , \\[3ex]\dis  
\dotfill  \  \\[3ex]\dis  
\frac{\Delta H}{\Delta y^{2m}} = \frac{H(y_{n+1}^1, y_{n+1}^2, \ldots, y_{n+1}^{2m}) - H (y_{n+1}^1, y_{n+1}^2, \ldots, y_n^{2m}) }{y_{n+1}^{2m} - y_n^{2m}} .
\ea \ee
where, to fix our attention, $\pmb{y}$ is defined by \rf{y=xp}. In fact, we may identify with $\pmb{y}$ any permutation of $2 m$ components $x^k, p^j$. Thus we have $(2 m)!$ discrete gradients of this type (in particular cases some of them may be identical). 

Having any discrete gradient we can easily obtain the related symmetric discrete gradient
\be  \label{grad-sym}
 {\bar \nabla}_s H (\pmb{y}_n, \pmb{y}_{n+1}) = \frac{1}{2} \left( {\bar \nabla} H (\pmb{y}_n, \pmb{y}_{n+1}) + {\bar \nabla} H (\pmb{y}_{n+1}, \pmb{y}_n) \right) \ .
\ee
One can easily verify that ${\bar \nabla}_s H$ satisfies conditions \rf{disgrad-def} provided that they are satisfied by $\bar \nabla H$.

\subsection{Linearization of discrete gradients} 

In order to construct locally exact modifications we need to linearize  discrete gradients defined by \rf{grad-incre} and \rf{grad-sym}. 

\begin{lem} \label{lem-inc}
Linearization of the coordinate increment discrete gradient \rf{grad-incre} around $\pmb{\bar y}$ yields
\be  \label{Tay-AB} 
  \bar \nabla H (\pmb{y}_n, \pmb{y}_{n+1}) \approx H_{\pmb{y}} + \frac{1}{2} \left( A \pmb{\nu}_{n+1} + B \pmb{\nu}_n \right) \ , 
\ee
where we denoted $\pmb{\nu}_n = \pmb{y}_n - \pmb{\bar y}$, and 
\be \ba{c} \dis  \label{AB}
A = \left(  \ba{lllll} \frac{1}{2} H_{y^1 y^1} & 0 & \ldots & 0 & 0 \\ H_{y^2 y^1} & \frac{1}{2} H_{y^2 y^2} & \ldots  & 0 & 0 \\ \multicolumn{5}{c}{\dotfill} \\ H_{y^{2 m-1} y^1} & H_{y^{2m-1} y^2} & \ldots & \frac{1}{2} H_{y^{2m-1} y^{2m-1}} & 0 \\ H_{y^{2m} y^1} & H_{y^{2m} y^2} & \ldots & H_{y^{2m} y^{2m-1} } & \frac{1}{2} H_{p^{2m} p^{2m}}   \ea              \right) \ , 
\\[9ex] \dis
B = \left(  \ba{lllll} \frac{1}{2} H_{y^1 y^1} & H_{y^1 y^2} & \ldots & H_{y^1 y^{2m-1}} & H_{y^1 y^{2m}} \\ 0 & \frac{1}{2} H_{y^2 y^2} & \ldots  & H_{y^2 y^{2m-1}} & H_{y^2 y^{2m}} \\ \multicolumn{5}{c}{\dotfill} \\ 0 & 0 & \ldots & \frac{1}{2} H_{y^{2m-1} y^{2m-1}} & H_{y^{2m-1} y^{2m}} \\ 0 & 0 & \ldots & 0 & \frac{1}{2} H_{y^{2m} y^{2m}}   \ea              \right) \ . 
\ea \ee
\end{lem}

\begin{Proof}
We  denote $\pmb{\hat y}_n^{j}= (y_{n+1}^1,\ldots,y_{n+1}^{j},y_n^{j+1},\ldots, y_n^{2 m})^T$. In particular, $\pmb{\hat y}_n^{0} = \pmb{y}_n$ and $\pmb{\hat y}_n^{2m} = \pmb{y}_{n+1}$. Expanding $H (\pmb{\hat y}_n^{j})$ around $\pmb{\hat y}_n^{j-1}$, we get
\be
 H (\pmb{\hat y}_n^{j}) = H (\pmb{\hat y}_n^{j-1}) + H_{y^j} (\pmb{\hat y}_n^{j-1}) (y_{n+1}^j - y_n^j) + \frac{1}{2} H_{y^j y^j} (\pmb{\hat y}_n^{j-1}) (y_{n+1}^j - y_n^j)^2 + \ldots  
\ee  
Hence
\be  \label{expa1}
\frac{\Delta H}{\Delta y^j} = \frac{ H (\pmb{\hat y}_n^{j}) - H (\pmb{\hat y}_n^{j-1})}{y_{n+1}^j - y_n^j} = 
H_{y^j} (\pmb{\hat y}_n^{j-1}) + \frac{1}{2} H_{y^j y^j} (\pmb{\hat y}_n^{j-1}) (y_{n+1}^j - y_n^j) + \ldots  
\ee
Then, from the definition of $\pmb{\nu}_n$, we have
\be \ba{l} \dis
y_{n+1}^j - y_n^j = \nu_{n+1}^j - \nu_n^j \ , \\[2ex]\dis
\pmb{\hat y}_n^{j-1} = \pmb{\bar y} + \left(  \nu_{n+1}^1,\ldots,\nu_{n+1}^{j-1},\nu_n^j,\ldots,\nu_n^{2 m} \right) \ ,
\ea \ee 
and, taking it into account, we rewrite \rf{expa1} as 
\be
\frac{\Delta H}{\Delta y^j} = H_{y^j} (\pmb{\bar y}) + \sum_{k=1}^{j-1} H_{y^j y^k} (\pmb{\bar y}) \ \nu_{n+1}^k + \sum_{k=j}^{2 m}  H_{y^j y^k} (\pmb{\bar y}) \ \nu_{n}^k + \frac{1}{2} H_{y^j y^j} (\pmb{\bar y}) (\nu_{n+1}^j - \nu_n^j) + \ldots 
\ee
which is equivalent to \rf{Tay-AB}, \rf{AB}. 
\end{Proof}

\begin{lem} \label{lem-sym}
Linearization of the symmetric discrete gradient yields
\be
\bar \nabla_s H (\pmb{y}_n, \pmb{y}_{n+1}) \approx  H_{\pmb{y}} + \frac{1}{2} H_{\pmb{y} \pmb{y}} \left( \pmb{\nu}_n + \pmb{\nu}_{n+1} \right) \ , 
\ee
where $H_{\pmb{y} \pmb{y}} $ is the Hessian matrix of $H$, evaluated at $\pmb{\bar y}$.  
\end{lem}

\begin{Proof}
We observe that  $A + B = H_{\pmb{y} \pmb{y}}$ 
and then we use \rf{grad-sym} and \rf{Tay-AB}.
\end{Proof}

\subsection{Conservative properties of modified discrete gradients}

In the one-dimensional case the corresponding locally exact modification is clearly energy-preserving, compare section~\ref{sec-grad1}. In the general case, conservative properties are less obvious. In this section we present  several useful results.

\begin{lem}  \label{lem-sigma}
We assume that a $2 m \times 2 m$ matrix $\Lambda$ (depending on $h$ and, possibly, on other variables) is skew-symmetric (i.e., $\Lambda^T = - \Lambda$) and 
\be  \label{AS} 
 \lim_{h \rightarrow 0} \frac{\Lambda}{h} = S \ , \qquad S = \mm 0 & 1 \\ - 1 & 0 \ema \ . 
\ee
Then, the numerical scheme  
\be  \label{num-A}
  \pmb{y}_{n+1} - \pmb{y}_n = \Lambda {\bar \nabla} H \ ,
\ee
where $\pmb{y}_n$ is defined by \rf{y=xp} and ${\bar \nabla} H$ satisfies \rf{disgrad-def}, is a consistent integrator for \rf{multiham}  preserving the energy integral up to round-off error. 
\end{lem}

\begin{Proof}
The consistency follows immediately form \rf{AS}. The energy preservation can be shown in the standard way. Using the standard scalar product in $\R^{2 m}$, we multiply both sides of \rf{num-A} by ${\bar \nabla} H$
\be
\scal{{\bar \nabla} H}{ \pmb{y}_{n+1} - \pmb{y}_n } = \scal{{\bar \nabla} H}{ \Lambda {\bar \nabla} H } \ .
\ee
By virtue of \rf{disgrad-def} the left-hand side equals $H (\pmb{y}_{n+1}) - H (\pmb{y}_n)$. The right hand side vanishes due to the skew symmetry of $\Lambda$. Hence $H (\pmb{y}_{n+1}) = H (\pmb{y}_n)$. 
\end{Proof}

\begin{lem}  \label{lem-zachow}
We assume that $2 m \times 2 m$ matrix \ $\pmb{\theta}$ \ is of the following form
\be  \label{theta}
  \pmb{\theta} = \mm \pmb{\delta} &  - \pmb{\sigma} \\ \pmb{\rho} & \pmb{\delta}^T \ema \ , \quad 
\pmb{\rho}^T = - \pmb{\rho} \ , \quad \pmb{\sigma}^T = - \pmb{\sigma} \ , \quad \lim_{h\rightarrow 0} \frac{\pmb{\theta}}{h} = 1 
\ee 
(where $\pmb{\delta}$, $\pmb{\sigma}$, $\pmb{\rho}$ are $m \times m$ matrices) and \ ${\bar \nabla} H $ \ is (any) discrete gradient. Then, the numerical scheme given by
\be  \label{mod-locex}
   \pmb{y}_{n+1} - \pmb{y}_n  = \pmb{\theta} S {\bar \nabla} H  
\ee
preserves the energy integral exactly, i.e., $H (\pmb{x}_{n+1}, \pmb{y}_{n+1}) = H (\pmb{x}_n, \pmb{y}_n)$.  
\end{lem}

\begin{Proof} We observe that \ $\pmb{\theta} S$ \ is skew-symmetric, and then we use Lemma~\ref{lem-sigma}.  
\end{Proof}

We easily see that \rf{mod-locex} reduces to the standard discrete gradient scheme when we take \ $\pmb{\theta} = h_n$ (i.e., \ $\pmb{\theta}$ is proportional to the unit matrix).

\begin{lem}  \label{lem-anal}
If \ $\pmb{\theta}$ is of the form \rf{theta} and $z \mapsto f (z)$ is any analytic function, then $f (\pmb{\theta})$ is also of the form \rf{theta}.    
\end{lem}

\begin{Proof}
First, we will show that the conditions \rf{theta}  are  equivalent to
\be  \label{theta1}
   \pmb{\theta}^T = S^{-1} \pmb{\theta} S \ . 
\ee
Indeed, assuming  a general form of $\pmb{\theta}$, e.g., $\pmb{\theta} = \mm \pmb{\delta} & - \pmb{\sigma} \\ \pmb{\rho} & \pmb{\gamma} \ema$ we see that  
the constraint \rf{theta1} is equivalent to $\pmb{\gamma} = \pmb{\delta}^T$, $\pmb{\rho}^T = - \pmb{\rho}$ and $\pmb{\sigma}^T = - \pmb{\sigma}$. Then the proof is straightforward. Assuming  $\dis f (z) = \sum_{k=1}^\infty a_n z^n$, we obtain
\be
 \left( \sum_{k=1}^\infty a_k \pmb{\theta}^k \right)^T = \sum_{k=1}^\infty a_k \left( S^{-1} \pmb{\theta} S \right)^k = S^{-1} \left( \sum_{k=1}^\infty a_k \pmb{\theta}^k \right) S \ , 
\ee
i.e., $ f (\pmb{\theta})^T = S^{-1} f (\pmb{\theta}) S$.  
\end{Proof}

\begin{cor}  \label{cor-pres}
If \ $\pmb{\theta}^T = S^{-1} \pmb{\theta} S$ \ and $f$ is an analytic function, then the scheme $\pmb{y}_{n+1} - \pmb{y}_n  = f (\pmb{\theta}) S {\bar \nabla} H $ preserves exactly the energy integral $H$.
\end{cor}

\subsection{Locally exact symmetric discrete gradient scheme}

We begin with the symmetric case because section~\ref{sec-locex-one} suggests that the coordinate increment discrete gradient case is more difficult.  In the symmetric case a locally exact modification is derived similarly as in the one-dimensional case. 

\begin{prop}  \label{prop-sym1-locex}
The following modification of the symmetric discrete gradient scheme  is locally exact at \ $\pmb{\bar y}$:  
\be  \label{sym1-locex} 
\pmb{y}_{n+1} - \pmb{y}_n = \pmb{\theta}_n S {\bar \nabla}_s H \ , 
\ee
where 
\be  \label{theta-n}
 \pmb{\theta}_n = 2 (F')^{-1} \tanh \frac{h_n F'}{2} \ , 
\ee
and $F'$, given by \rf{F-F'}, is evaluated at $\pmb{y} = \pmb{\bar y}$. 
\end{prop}

\begin{Proof}  We are going to derive \rf{theta-n}, assuming that $\pmb{\theta}_n$ depends on $\pmb{\bar y}$ and $h$. By virtue of Lemma~\ref{lem-sym} the linearization of \rf{sym1-locex} is given by
\be \label{lin1} 
\pmb{\nu}_{n+1} - \pmb{\nu}_n = \pmb{\theta}_n S \left( H_{\pmb{y}} + \frac{1}{2} H_{\pmb{y} \pmb{y}} (\pmb{\nu}_{n+1} + \pmb{\nu}_n) \right) \ .
\ee
Taking into account \rf{F-F'} we transform \rf{lin1} into 
\be \label{nu-dis} 
 \left( 1 - \frac{1}{2} \pmb{\theta}_n F' \right) \pmb{\nu}_{n+1} =  \left( 1 + \frac{1}{2} \pmb{\theta}_n F' \right) \pmb{\nu}_n + \pmb{\theta}_n F \ .
\ee
The scheme \rf{sym1-locex} is locally exact iff \rf{nu-dis} coincides with \rf{lin-hex}. Therefore, we require that 
\be   \label{two-eq1}
\left( 1 - \frac{1}{2} \pmb{\theta}_n F' \right) e^{h_n F'} =  1 + \frac{1}{2} \pmb{\theta}_n F'  \ , 
\ee
\be  \label{two-eq2}
\left( 1 - \frac{1}{2} \pmb{\theta}_n F' \right) \left( e^{h_n F'} - 1 \right) (F')^{-1} F = \pmb{\theta}_n F \ .
\ee
From equation \rf{two-eq1} we can compute $\pmb{\theta}_n$ which yields \rf{theta-n}. Equation \rf{two-eq2} is automatically satisfied provided that \rf{two-eq1} holds. 
\end{Proof}

\begin{prop}  \label{prop-sym1-pres}
The numerical scheme \rf{sym1-locex} with $\pmb{\theta}_n$ given by \rf{theta-n} is energy-preserving. 
\end{prop}

\begin{Proof} We have 
$F' = S H_{\pmb{y} \pmb{y}}$. Therefore, $(F')^T = - H_{\pmb{y} \pmb{y}} S = - S^{-1} F' S$, and, as a consequence 
\be  \label{F'2}
((F')^2 )^T = S^{-1} (F')^2 S \ . 
\ee
It means that $(F')^2$ has the form \rf{theta}. The formula \rf{theta-n} expresses $\pmb{\theta}_n$ as an analytic function of $(F')^2$.  Finally, we use Lemma~\ref{lem-anal}. 
\end{Proof}

In the one-dimensional case $(F')^2$ is proportional to the unit matrix which essentially simplifies arguments presented in this section.

\subsection{Locally exact coordinate increment discrete gradient scheme}

The symmetric form of the discrete gradient leads to a simple form of the locally exact modification. It turns out, however, that starting from the simplest form of the discrete gradient, namely coordinate increment discrete gradient \cite{IA}, we also are able to derive the corresponding locally exact modification.

\begin{prop}
The following modification of the coordinate increment discrete gradient scheme is locally exact at \ $\pmb{\bar y}$:  
\be  \label{incre1-locex} 
\pmb{y}_{n+1} - \pmb{y}_n = \pmb{\theta}_n S {\bar \nabla} H \ , 
\ee
where ${\bar \nabla} H$ is given by \rf{grad-incre}, 
\be  \label{theta2-n}
 \theta_n = 2 \left( S R + F' \coth \frac{h_n F'}{2} \right)^{-1} \ , 
\ee
$F'$ is given by \rf{F-F'} (i.e., $F' = S H_{\pmb{y} \pmb{y}}$), and, finally $R = A - B$, i.e., 
\be
 R = \left(  \ba{ccccc} 0  & - H_{y^1 y^2} & \ldots & - H_{y^1 y^{2m -1}} & - H_{y^1 y^{2 m}} \\ H_{y^2 y^1} & 0 & \ldots  & - H_{y^2 y^{2m -1}} & - H_{y^2 y^{2 m}} \\ \multicolumn{5}{c}{\dotfill} \\ H_{y^{2 m-1} y^1} & H_{y^{2m-1} y^2} & \ldots & 0 & H_{y^{2m-1} y^{2m}}  \\ H_{y^{2m} y^1} & H_{y^{2m} y^2} & \ldots & H_{y^{2m} y^{2m-1} } & 0    \ea    \right)  \ .     
\ee
$F'$ and $R$ are evaluated at \ $\pmb{y} = \pmb{\bar y}$. 
\end{prop}

\begin{Proof} We are going to derive \rf{theta2-n}, assuming that $\pmb{\theta}_n$ depends on $\pmb{\bar y}$ and $h$. By virtue of Lemma~\ref{lem-inc} the linearization of \rf{incre1-locex} is given by
\be  \label{lin2}
 \pmb{\nu}_{n+1} - \pmb{\nu}_n = \pmb{\theta}_n S ( A \pmb{\nu}_{n+1} + B \pmb{\nu}_n ) + \pmb{\theta}_n  S H_{\pmb{y}} \ .
\ee
Hence, taking into account that $S H_{\pmb{y}} = F$,   
\be \label{nu-dis2} 
 \left( 1 -  \pmb{\theta}_n  S A \right) \pmb{\nu}_{n+1} =  \left( 1 +  \pmb{\theta}_n  S B \right) \pmb{\nu}_n + \pmb{\theta}_n F \ .
\ee
The scheme \rf{incre1-locex} is locally exact iff \rf{nu-dis2} coincides with \rf{lin-hex}. Therefore, we require that 
\be   \label{2-eq1}
\left( 1 - \pmb{\theta}_n  S A \right) e^{h_n F'} =  1 +  \pmb{\theta}_n  S B  \ , 
\ee
\be  \label{2-eq2}
\left( 1 -  \pmb{\theta}_n  S A \right) \left( e^{h_n F'} - 1 \right) (F')^{-1} F = \pmb{\theta}_n F \ .
\ee
Inserting \rf{2-eq1} into \rf{2-eq2} we get
\be
\pmb{\theta}_n  S ( B  + A ) (F')^{-1} F = \pmb{\theta}_n F \ ,
\ee
which is identically satisfied by virtue of 
 $A + B = H_{\pmb{y} \pmb{y}} = S^{-1} F'$,  compare \rf{AB}. The remaining equation, \rf{2-eq1}, defines $\pmb{\theta}_n$: 
\be  \label{defteta}
\pmb{\theta}_n  \left( S A   e^{h_n F'}  + S B \right) =  e^{h_n F'} - 1 \ .
\ee
In order to simplify \rf{defteta} we introduce $R = A - B$ ($R$ is antisymmetric because $B = A^T$). Then, taking into account $SA+SB= F'$, we get
\be
 S A  = \frac{1}{2}  F'  + \frac{1}{2} S R \ , \qquad S B = \frac{1}{2}  F'  - \frac{1}{2} S R \ . 
\ee
Substituting it into \rf{defteta} we complete the proof. 
\end{Proof}

\begin{prop}
The numerical scheme \rf{incre1-locex} with $\pmb{\theta}_n$ given by \rf{theta2-n} is energy-preserving, i.e., $H (\pmb{y}_{n+1}) = H (\pmb{y}_n)$. 
\end{prop}

\begin{Proof} We have 
\be 
  \pmb{\theta}_n^T = 2 \left( (S R)^T + \left( F' \coth \frac{h_n F'}{2} \right)^T \right)^{-1} = S^{-1}  \pmb{\theta}_n S \ ,
\ee
because
\be
   (S R)^T = (- R)(-S) = S^{-1} \left( S R \right) S 
\ee
and,  by virtue of Lemma~\ref{lem-anal}, 
\be
\left( F' \coth \frac{h_n F'}{2} \right)^T = S^{-1} \left( F' \coth \frac{h_n F'}{2} \right) S \ , 
\ee 
where we took into account \rf{F'2}. Then, we use Corollary~\ref{cor-pres}. 
\end{Proof}

\subsection{Separable Hamiltonians. Multidimensional case}
\label{sec-multli-sep}

In the case of one degree of freedom the assumption $H_{\pmb{x} \pmb{p}} = 0$ simplifies final formulas and yields the same results for both considered gradient schemes, see section~\ref{sec-one-sep}. In the multidimensional  separable case, i.e., 
\be  \label{H-sep} 
H (\pmb{x}, \pmb{p}) = T (\pmb{p}) + V (\pmb{x}) \ , 
\ee
 a considerable simplification occurs only for the symmetric discrete gradient.

\begin{prop} \label{prop-TV}
The numerical scheme 
\be \ba{l}
\pmb{\delta}_n^{-1} \left( {\pmb x}_{n+1} - {\pmb x}_n \right) = {\bar \nabla}_s T (\pmb{p}_n, \pmb{p}_{n+1})    \\[2ex] \dis
(\pmb{\delta}_n^T)^{-1} \left( {\pmb p}_{n+1} - {\pmb p}_n \right) = - {\bar \nabla}_s V (\pmb{x}_n, \pmb{x}_{n+1})    
\ea \ee
preserves exactly the energy integral (for any  $m\times m$ matrix $\pmb{\delta}_n$), i.e., $T (\pmb{p}_n) + V (\pmb{x}_n)$ does not depend on $n$. This scheme is locally exact for
\be
   \pmb{\delta}_n = 2 \Omega_n^{-1} \tan \frac{h_n \Omega_n}{2} \ , \qquad \Omega_n^2  = T_{\pmb{p} \pmb{p}} (\pmb{\bar x}, \pmb{\bar p}) V_{\pmb{x} \pmb{x}} (\pmb{\bar x}, \pmb{\bar p})  \ . 
\ee
\end{prop}

\begin{Proof} 
The first part of the Proposition follows directly from Lemma~\ref{lem-zachow} (in this case $\pmb{\sigma} = \pmb{\rho} = 0$, and $\pmb{\delta} = \pmb{\delta}_n$). The second part is a consequence of Proposition~\ref{prop-sym1-locex}. 
We have
  $F = (T_{\pmb{p}},  - V_{\pmb{x}})^T$ \ and  
\be
 F' = \mm  0 &  T_{\pmb{p} \pmb{p}} \\ 
- V_{\pmb{x} \pmb{x}} & 0 \ema  , \qquad (F')^2 = - \mm T_{\pmb{p} \pmb{p}} V_{\pmb{x} \pmb{x}} & 0 \\ 0 & V_{\pmb{x} \pmb{x}} T_{\pmb{p} \pmb{p}} \ema  ,
\ee
where all quantities are evaluated at $(\pmb{\bar x}, \pmb{\bar p})$. 
We denote $\Omega_n^2 = T_{\pmb{p} \pmb{p}} V_{\pmb{x} \pmb{x}}$, and then
\be
 (F')^2 = - \mm \Omega_n^2 & 0 \\ 0 & (\Omega_n^T)^2 \ema
\ee
and
\be
 2 (F')^{-1} \tanh \frac{h_n F'}{2} = \mm  h_n \tanc  \left( \frac{1}{2} h_n \Omega_n \right) & 0  \\ 0 & h_n \tanc  \left( \frac{1}{2} h_n \Omega_n^T \right)  \ema \ .
\ee
We complete the proof applying Proposition~\ref{prop-sym1-locex}, compare also 
Proposition~\ref{prop-sym1-pres}.
\end{Proof}

 Formulas of Corollary~\ref{cor-harm} are strikingly similar to those given in Proposition~\ref{prop-TV}.   
This is due to the fact that locally exact discretizations applied to linear systems yield exact integrators. 
Indeed, the harmonic oscillator \rf{vec1-osc}  is a special case of \rf{H-sep} for 
\be
T (\pmb{p}) = \frac{1}{2} \pmb{p}^2 \ , \qquad  V (\pmb{x}) = \frac{1}{2} \scal{\pmb{x}}{\Omega^2 \pmb{x}} - \scal{\pmb{x}}{\pmb{a}} 
\ee
Then $T_{\pmb{p} \pmb{p}} = 1$ and $V_{\pmb{x} \pmb{x}} = \Omega^2$. Hence  $\Omega^T = \Omega$ and $\pmb{\delta}_n^T = \pmb{\delta}_n$. Proposition~\ref{prop-harm} is an obvious consequence of Proposition~\ref{prop-TV}.

\section{Concluding remarks}

We presented a new construction of very accurate numerical schemes based on the notion of local exactness. This notion is known (although under different names) since almost fifty years. The original application, see \cite{Pope}, has been confined to the exact discretization of linearized equations, compare section~\ref{sec-ex-lin}. We obtain in this way one particular locally exact integrator, which has some advantages (e.g., it can serve as a good predictor, (\cite{CR-PRE}, section V.A) but lacks geometric properties and stability  of, for instance, locally exact discrete gradient methods \cite{CR-PRE}. 

Our approach has two new features. First, we modify known numerical schemes in a locally exact way. Any numerical scheme admits at least one (usually more) natural locally exact modification, see section~\ref{sec-pop}. Second, we try to preserve geometric properties of the original numerical scheme. This task is not trivial. In this paper we present one successful application: locally exact modifications of discrete gradient methods for canonical Hamilton equations. 
We are able to preserve exactly the energy integral and, in the same time, increase the accuracy by many orders of magnitude. Another advantage is a variable time step. Unlike symplectic methods (which work mostly for the constant time step) discrete gradient methods admit conservative modifications with variable time step. Therefore, one may easily  implement any variable step method in order to obtain further improvement. 

We devoted a lot of attention to the one-dimensional case, presenting it independently from the general, multidimensional case. The reason is not only pedagogical but also practical. One-dimensional case already proved to be successful \cite{CR-PRE,CR-BIT}. The proposed modification, although more expensive (but only by a dozen or so percents), turns out to be more accurate even by 8 orders of magnitude in comparison to the standard discrete gradient scheme. Therefore, in the case of one degree of freedom our modifications are very efficient. In multidimensional cases the relative cost of our algorithm is higher, but still we hope that our method will be of advantage. 
Modifications proposed in this paper contain exponentials of variable matrices. Similar time-consuming evaluations are characteristic for all exponential integrators and  in this context 
effective methods of computing matrix exponentials have been recently developed \cite{HoL1,NW}. 

We presented locally exact modifications from a unified theoretical perspective. There are many possible further developments. First of all, we plan to apply our approach to chosen multidimensional problems, testing the accuracy of locally exact schemes by numerical experiments. Then, we would like to extend the range of applications. 
Throughout this paper we assumed the autonomous case, ${\pmb {\dot x}} = F ({\pmb x})$. The extension on the non-autonomous case can be done along lines indicated already in the Pope's paper \cite{Pope}. A separate problem is to obtain in this case any  locally exact modification with geometric properties.
Another open problem is the construction of locally exact (or, at least, linearization-preserving) deformations of generalized discrete gradient algorithms preserving all first integrals (see \cite{MQR2}).  
Linearization-preserving integrators form an important subclass of locally exact schemes. It would be worthwhile to study linearization-preserving modifications of geometric numerical integrators, especially in those cases when locally exact are difficult or impossible to construct.

\end{document}